# The calculation of $p_n$ and $\pi(n)$

Simon Plouffe
June 03, 2020

Abstract

A new approach is presented for the calculation of $p_n$ and $\pi(n)$ which uses the Lambert W function. An approximation is first found and using a calculation technique it makes it possible to have an estimate of these two quantities more precise than those known from Cipolla and Riemann. The calculation of $p_n$ uses an approximation using the Lambert W function and an estimate based on a logarithmic least square curve (LLS) c (n). The formula is:

$$p_n + \pi(n) \approx -nW_{-1}\left(\frac{-e}{n}\right) c(n) \qquad 1$$

The results presented are empirical and apply up to $n \approx 1.358 \times 10^{16}$.

# Introduction

Today we know very large prime numbers like $2^{82589933} - 1$ but the rank of these numbers is unknown. The primes with their rank is known up to $10^{27}$ for $\pi(n)$ and $10^{24}$ for $p_n$ at some specific points and the complete list of primes with their rank is known only up to $10^{17}$.

In 2010, Dusart [16] proved that $\pi(n) \approx \frac{n}{\log(n)-1}$ if $n > 5393$. We will use this approximation to give one approximation of $p_n$ by inverting the formula.

$$\text{If } \pi(n) \approx \frac{n}{\ln(n)-1} \text{ then } p_n \approx -nW_{-1}\left(\frac{-e}{n}\right).$$

$W_{-1}(n)$ is the Lambert W function of order -1, $W_0$ or $W(n)$ is of order 0. We We take the branch of the W function which gives meaning to quantity $\frac{-e}{n}$.

That formula is quite precise, for $n = 10^{24}$, we have $p_n$ precise to 99.97 %.

By analyzing the rest of $p_n$ and $-nW_{-1}\left(\frac{-e}{n}\right)$, we rapidly find that it is close to $\frac{n}{W(n)}$, so $p_n \approx -nW_{-1}\left(\frac{-e}{n}\right) - \frac{n}{W(n)}$. But also that the quantities

$$\frac{n}{\ln(n)-1} \approx Li(n) \approx \pi(n) \approx \frac{n}{W(n)} \approx \frac{\rho_n}{2\pi}$$

are close too when $n \to \infty$. Here, $\rho_n$ is the imaginary part of the n'th non-trivial zero of Riemann's $\zeta$ and $Li(n)$ is the integral logarithm, the values are : .1843e23, .18434e23, .1844e23, .1948e23 et .1986e23 respectively.

Empirical data suggest that the choice of $\pi(n)$ is the best.

If $n = 10^{24}$ then

$$p_{10^{24}} + \pi(10^{24}) \approx -10^{24}W_{-1}\left(\frac{-e}{10^{24}}\right) \tag{2}$$

The approximation is valid up to 0.99999401 or 99.999401 % of the real value.
A naive question is then, couldn't we be more precise? Generally speaking, is there a way to have an exact formula? If we can be more precise then how much? What about the rest? What is his nature exactly?

## An explicit formula: Prime numbers in geometric progression

Right now, the question of whether an exact formula exists for prime numbers is yes and no both. There are many formulas that give them but they are either impractical or limited see the table in appendix.

For example, here is an explicit formula which gives some primes in geometric progression.

$$f(n) = \{c^n\} \tag{3}$$

Here { } is the nearest integer function. The smallest example of c is : if c = 2.553854696... then $f(n)$ is 3, 7, 17, 43, 109, 277, 709. The sequence contains only 7 terms. But can we go further? Using the simulated annealing method and Monte Carlo we can go much further.

| Constant c such that $\{c^n\}$ is prime | Values | Number of primes generated |
|---|---|---|
| 2.553854696... | 3, 7, 17, 43, 109, 277, 709 | 7 |
| 593.46526943871... | n=2..48 | 47 |
| 2027.167168476491219434395 | n=1..97 | 97 |
| 29983.279826631136... | n=1..422 | 422 |
| 55237.075042967647154331247815286617... | n=2..633 | 632 |

From there we can conjecture that if c is large enough, the sequence with formula (3) can generate an arbitrary number of primes explicitly, see [23].

All we have for the expression of $p_n$ or $\pi(n)$ for large values of n are just approximations. The best known approximation is that of Riemann, the function is called Riemann-R or that of Gram found in 1884. The only way that has been found is to use partially Gram's formula followed by a sophisticated Erathosthenes screen. This is why that the value of $p_n$ or $\pi(n)$ is only known up to $10^{24}$ and $10^{27}$ respectively. The last calculation of $\pi(10^{27})$ took several months and a significant amount of computer time (23 CPU years) in 2018.

# First approximation

We will confine ourselves for the moment to the calculation of $p_n$ since $\pi(n)$ can be calculated by inversion.
The classic formula for $p_n$ is $p_n \approx n \ln(n)$ or better yet the one found by Cipolla in 1902 states that

$$p_n \approx n\left(\ln(n) + \ln(\ln(n)) - 1 + \frac{\ln(\ln(n)) - 2}{\ln(n)} - \frac{\ln(\ln(n))^2 - 6\ln(\ln(n)) + 11}{2\ln(n)^2} + \cdots \right). \quad (4)$$

The calculation was taken further in 1994 by B. Salvy who extracted a procedure for pushing it further the approximation.
What is remarkable is the similarity with the asymptotic development of $W(n)$.

$$W(n) \approx L_1 - L_2 + \frac{L_2}{L_1} + \frac{L_2(-2 + L_2)}{2L_1^2} + \frac{L_2(6 - 9L_2 + 2L_2^2)}{6L_1^3} + \frac{L_2(-12 + 36L_2 - 22L_2^2 + 36L_2^3)}{12L_1^4} + \cdots$$

$L_1 = \ln(n)$ et $L_2 = \ln(\ln(n))$.

The calculation of $p_n$ with this formula (4) gives 12 exact decimal places out of the 26 that count $p_{10^{24}}$. But does not allow to go further, even with 64 terms in the asymptotic development. The formula of $Gram^{<-1>}$ is much more precise. $Gram^{<-1>}$ is the functional inverse of the Gram series [24].

# A better approximation

A summary analysis indicates that the rest after the first term of (1) is of type logarithmic : $a + b\,ln(n)$ where a and b to be determined. An idea is then to calculate the logarithmic least squares curve (LLS) passing through a chosen number of points on a table of values of primes.

We can also notice that by taking only one term for the approximation of $-nW_{-1}(\frac{-e}{n})$, this form is equivalent to several terms of Cipolla's development. If we take the 2 terms it will be even more precise. In other words, given the nature of the asymptotic development of $W(n)$, each term is equivalent to several terms in the development of Cipolla. We assume here that the remainder after the 2 terms is a logarithmic curve and that once calculated it will stick to reality.

The question then arises of what is the nature of what remains? In fact, we do not know. The best known for $\pi(n)$ which is theoretically valid is $Li(n)$. Riemann proposed a 2nd formula that is much better at first but was invalidated by Littlewood in 1914. This 2nd formula, called Riemann R or equivalent, the Gram series is

$$\pi(n) = \sum_{k=1}^{\infty} \frac{\ln(x)^k}{k\,k!\,\zeta(k+1)}$$

Numerically, the approximation of $\pi(n)$ by Riemann's R or the Gram series is excellent in addition to converging quickly. But Littlewood showed that after $10^9$, the approximation is less precise. As for the function $Li(n)$, it behaves better at much larger scales, the first crossover of $Li(n)$ and $\pi(n)$ was evaluated at $10^{316}$ approximately, that is $Li(n) - \pi(n) = 0$ in the neiboughoud of $1.397 \times 10^{316}$.

The best approximation that has been found empirically is:

$$p_n = -nW_{-1}\left(\frac{-e}{n}\right)c(n) - \pi(n)$$

where $c(n)$ is of type $a + b\,log(n)$.

Here is a table of values that were considered for the calculation of $p_n$ from $10^2$ to $10^{17}$

| Value of the step | Number of values | Range |
|---|---|---|
| $10^2$ | 27117419 | 2711741900 |
| $10^3$ | 32082085 | 32082085000 |
| $10^4$ | 45020269 | 450202690000 |
| $10^5$ | 16038989 | $1.603 \times 10^{12}$ |
| $10^6$ | 4046531 | $4.046531 \times 10^{12}$ |
| $10^7$ | 5011691 | $5.011691 \times 10^{13}$ |
| $10^8$ | 454060 | $4.54060 \times 10^{13}$ |
| $10^9$ | 2200000 | $2.2 \times 10^{15}$ |
| $10^{10}$ | 1358121 | $1.358121 \times 10^{16}$ |
| $10^{11}$ | 135812 | $1.358121 \times 10^{16}$ |
| $10^{12}$ | 54974 | $5.4974 \times 10^{16}$ |
| $10^{13}$ | 12317 | $1.2317 \times 10^{17}$ |
| $10^{14}$ | 2162 | $2.162 \times 10^{17}$ |

We can calculate the points on the LLS curve by solving the following eqaution for x, taking n in one of the tables above.

$$-nW_{-1}\left(\frac{-e}{n}\right)x - \pi(n) + p_n = 0$$

by using the bisection method. The values are between 0.9 and 1. We then calculate the logarithmic least squares curve. The coefficient $r^2$ will indicate if the curve is right. The calculation of coefficients a and b is done according to the formula:

$$b = \frac{n\sum_{i=1}^{n}(y_i \ln x_i) - \sum_{i=1}^{n} y_i \sum_{i=1}^{n} \ln x_i}{n\sum_{i=1}^{n}(\ln x_i)^2 - (\sum_{i=1}^{n} \ln x_i)^2}$$

$$a = \frac{\sum_{i=1}^{n} y_i - b\sum_{i=1}^{n}(\ln x_i)}{n}$$

Recall, the coefficient $r^2$ indicates whether the experimental data stick to the curve. If is near 1 or 0, the curve follows a straight line very closely. The logarithmic least-squares curve is simply the log of values that are aligned on a line.

The curve thus found for $n \leq 1.358121 \times 10^{16}$ is

$$c(n) = 1.0000314775792421150615325693061 \\ - 0.00000051483940138413674623044640769440 \, ln(n)$$

Once you have found this LLS curve, you can significantly increase the accuracy if you use a trick. We take the chosen interval, $n \cdot 10^{10}$, $n = 1..1358121$. We then separate into 1358 slices of 10,000 values and for each value we modify the curve $c(n)$ with the formula

$$c(n) \to c(n) \, s^k$$

And $s = (1 - 10^{-10})$, then it suffices to find for which k, the curve admits a minimum error if we compare to the true values of $p_n$. An experiment was conducted with 13581 intervals to see if the mean value of the deviation decreased: this is the case. The reasonable limit found here is from 1,358 slices of 10,000 values.
In the appendix you can consult the program written in Maple language which performs the operation.

# Calculation of $\pi(n)$

To calculate $\pi(n)$, it suffices to isolate $\pi(n)$ in the equation of $G(n)$, normally trivial but in practical it does not work. Indeed $\pi(n)$ is smaller than $p_n$ in size. The formula remains valid except that the coefficients change slightly. The calculation of can be done with formula 1, with the primecount program which gives it directly. We can also have the values of Gram(n) and Gram inverted. If instead of $\pi(n)$ in formula (1) we rather take the term $\frac{n}{W(n)}$, it's a little less precise but still allows good accuracy. For example, for the interval from 1 to 1000000 the following program calculates $p_n$ very precisely.

```
G:=proc(n) # calculation of p(n) up to 1 million (first million of primes).
local ll,lk,s,s2,ss,kk;
  ll := [82, 16, -14, 4, -31, -10, -31, -32, -1, -44, -17, -38, -8, 7,
    -35, -41, -38, -3, 6, -14, -27, -12, -5, -51, -40, 17, -7, -17, -16,
    -14, 13, -7, -5, -1, -26, -29, -27, -31, -9, 8, 16, 4, 9, 0, 20, 11,
    7, -15, -23, -17, -10, -2, 2, -5, 8, 7, 9, 3, -12, -11, 4, 5, -3, 9,
    -1, 7, 24, 25, 33, 20, 15, 11, 9, 3, 9, 15, 3, 3, 1, 4, 5, -9, -1,
    -12, -4, 14, 16, 17, 28, 18, 12, 21, 24, 10, 14, 16, 15, 24, 26, 36,30];
  lk := floor(1/10000*n) + 1;
  s  := .3238679016803340+.4042167153029803e-1*ln(kk):
  s2 := subs(kk = n, s);
  ss := s2*0.999^ll[lk];
  round(abs(-n*W(-1, -exp(1.0)/n)) - ss*fab(n/W(n)));
end;
```

By inverting $\pi(n)$ to find , over the interval 1 to 1350000x10e16. We have the following program. It has the advantage of not having to calculate or having a table of the values of the first, but it is less precise.

```
f:=proc(k) local n,a,bas,haut,ll,lk,s,ss; # calcul de pi(n) jusqu'a 10^16
haut:=evalf(2*k/log(k)):
bas:=haut/8;

ll:=[3122, 3186, 3222, 3243, 3251, 3270, 3272, 3282, 3289, 3291, 3289, 3297, 3312,
3305, 3300, 3313, 3308, 3311, 3321, 3318, 3323, 3322, 3326, 3322, 3319, 3322,
3328, 3328, 3327, 3335, 3334, 3336, 3336, 3337, 3335, 3335, 3331, 3334, 3337,
3334, 3341, 3343, 3345, 3346, 3343, 3340, 3341, 3342, 3348, 3348, 3347, 3349,
3349, 3345, 3349, 3350, 3348, 3352, 3351, 3346, 3344, 3346, 3344, 3348, 3348,
3349, 3354, 3355, 3357, 3355, 3354, 3355, 3358, 3359, 3358, 3358, 3358, 3357,
3357, 3358, 3356, 3355, 3357, 3358, 3359, 3358, 3355, 3355, 3360, 3357, 3354,
3358, 3358, 3362, 3361, 3360, 3360, 3361, 3360, 3360, 3362, 3361, 3362, 3363,
3363, 3362, 3364, 3363, 3363, 3361, 3360, 3361, 3364, 3366, 3366, 3366, 3366,
3368, 3369, 3369, 3368, 3366, 3367, 3368, 3367, 3368, 3368, 3368, 3368, 3366,
3366, 3365, 3365, 3364, 3365]:

lk:=floor(k/100000000000000)+1:

s:=.882819461483173314372633+.855943969749036445417381e-3*ln(n):
ss:=s*(0.999999)^ll[lk]:

a:=abs(evalf(-n*W(-1,-exp(1.0)/n)))-evalf(ss)*fab(n/W(n))-k;
fsolve(a,n,n=bas..haut,fulldigits);
end;
```

# Appencix (Tables and Programs)

## Calculation of $p_n$

## Compison with the range 10000000000 ($10^{10}$) ... $1.352 \times 10^{16}$

| Formula for $p_n$ | Formula of Gram inverted | G(n) | Formula of Cipolla-Salvy |
|---|---|---|---|
| Minimal gap | 57 | 13 | 640495 |
| Maximal gap | 117539110 | 412614395 | 1103 millions |
| Average gap | 79.23 millions | 18.81 millions | 510 millions |

Conclusion : The formula with LambertW function is 4.21 times more precise than that of Gram series (inverted).

# Maple program for the calculation of $p_n$, n ≤ 1.356 × $10^{16}$

################################################################
Digits:=32:
**G:=proc(n, piofn)**
local cn, ll, s, ss, lk, s2;
  ll := [800, 369, 97, -19, -69, -133, -151, -138, -195, -200, -161, -182, -210, -177, -212,
    -190, -195, -138, -169, -199, -205, -199, -177, -175, -157, -150, -165, -156, -163,
    -162, -136, -158, -169, -175, -152, -150, -121, -134, -143, -134, -120, -125, -107,
    -115, -130, -90, -108, -125, -132, -134, -135, -134, -125, -119, -103, -88, -62, -88,
    -85, -89, -93, -89, -80, -77, -77, -74, -83, -70, -91, -89, -82, -75, -80, -78, -83,
    -82, -68, -61, -56, -47, -50, -63, -65, -74, -77, -72, -63, -67, -69, -72, -64, -41,
    -31, -29, -28, -28, -37, -42, -36, -40, -27, -26, -16, -22, -31, -37, -41, -38, -47,
    -40, -39, -39, -41, -42, -38, -37, -33, -39, -42, -31, -32, -39, -30, -25, -22, -11,
    -18, -18, -22, -18, -25, -30, -25, -24, -18, -15, -9, -6, -6, -8, -5, -1, 1, -14, -8,
    -9, -4, 0, -4, 4, 5, 4, 3, -2, -4, -10, 0, 1, 0, 4, -3, -8, -11, -13, -3, -1, 4, 1, -5,
    0, 5, 6, 3, -2, -2, 3, 4, -2, 2, 2, 2, 2, -2, -1, -4, 3, 1, -1, 3, 5, 16, 14, 12, 18,
    14, 13, 15, 15, 17, 17, 13, 12, 3, 3, 5, 8, 8, 9, 11, 14, 19, 22, 20, 17, 14, 16, 16,
    19, 17, 23, 28, 29, 20, 22, 16, 19, 22, 21, 17, 19, 26, 28, 28, 28, 26, 25, 20, 19, 20,
    20, 24, 24, 31, 35, 31, 28, 33, 31, 31, 31, 28, 31, 30, 30, 41, 39, 41, 34, 32, 24, 27,
    26, 29, 30, 29, 23, 21, 19, 20, 14, 17, 17, 11, 13, 19, 20, 19, 19, 16, 18, 16, 17, 22,
    28, 31, 35, 35, 37, 34, 31, 30, 33, 25, 26, 26, 23, 21, 22, 21, 24, 23, 19, 22, 23, 24,
    27, 24, 23, 23, 25, 25, 21, 27, 28, 34, 34, 33, 36, 35, 36, 37, 37, 34, 32, 33, 38, 47,
    48, 46, 45, 47, 44, 44, 50, 43, 40, 44, 46, 44, 39, 41, 42, 43, 43, 45, 44, 42, 44, 40,
    41, 42, 42, 42, 40, 40, 40, 36, 35, 37, 36, 37, 37, 37, 32, 35, 37, 37, 38, 31, 31,
    31, 28, 25, 28, 23, 25, 25, 23, 26, 29, 30, 29, 28, 26, 27, 30, 35, 35, 35, 35, 34, 38,
    37, 36, 35, 37, 39, 40, 35, 39, 42, 38, 37, 40, 42, 42, 43, 42, 40, 39, 39, 40, 40, 38,
    37, 40, 41, 41, 36, 34, 38, 38, 34, 34, 33, 36, 35, 34, 34, 38, 38, 37, 36, 39, 37, 37,
    32, 33, 31, 31, 27, 28, 23, 27, 25, 25, 30, 30, 31, 30, 31, 30, 35, 34, 32, 34, 37, 37,
    37, 40, 42, 38, 37, 36, 36, 36, 39, 37, 37, 38, 37, 34, 31, 32, 32, 30, 29, 29, 29, 30,
    25, 26, 27, 28, 24, 22, 22, 26, 29, 33, 35, 36, 35, 34, 35, 36, 36, 35, 33, 30, 34, 36,
    35, 37, 34, 38, 36, 36, 35, 36, 36, 37, 38, 38, 36, 36, 37, 39, 37, 38, 36, 37, 35, 36,
    38, 36, 35, 37, 39, 38, 35, 34, 34, 34, 35, 36, 33, 35, 34, 36, 39, 44, 44, 40, 37, 38,
    39, 35, 35, 34, 34, 37, 38, 37, 36, 36, 33, 30, 31, 29, 30, 33, 33, 31, 30, 31, 32, 33,
    31, 33, 33, 32, 31, 31, 29, 29, 29, 34, 33, 33, 31, 31, 32, 31, 31, 32, 33, 34, 33, 32,
    33, 35, 35, 31, 30, 32, 32, 32, 32, 33, 33, 32, 31, 30, 31, 32, 31, 27, 28, 28, 29, 31,
    32, 33, 30, 29, 29, 29, 27, 27, 28, 29, 29, 28, 28, 27, 29, 27, 27, 28, 26, 28, 27, 25,
    25, 26, 27, 27, 25, 24, 23, 20, 20, 20, 21, 23, 21, 21, 21, 22, 22, 23, 24, 22, 22, 22,
    20, 23, 20, 19, 19, 19, 20, 21, 23, 23, 23, 22, 18, 17, 16, 15, 12, 12, 12, 15, 13, 15,
    12, 11, 8, 13, 14, 15, 16, 18, 18, 17, 16, 19, 21, 21, 21, 22, 23, 22, 22, 23, 24, 24,
    22, 23, 23, 23, 22, 21, 20, 20, 20, 19, 18, 16, 14, 15, 18, 19, 17, 15, 16, 16, 19, 19,
    19, 19, 20, 20, 16, 16, 18, 19, 19, 17, 18, 17, 17, 16, 17, 17, 18, 19, 20, 20, 22, 21,
    21, 21, 20, 18, 18, 18, 19, 19, 19, 21, 20, 20, 21, 21, 23, 23, 22, 20, 17, 16, 17, 15,
    14, 14, 14, 16, 17, 18, 17, 17, 17, 18, 17, 18, 18, 17, 17, 18, 17, 17, 18, 18, 19, 17,
    17, 16, 15, 15, 16, 14, 14, 15, 15, 15, 14, 14, 14, 14, 15, 15, 13, 16, 16, 14, 15, 15,
    14, 13, 12, 11, 11, 12, 10, 9, 10, 10, 10, 9, 9, 10, 11, 11, 11, 12, 11, 11, 8, 8, 9,
    10, 10, 11, 10, 10, 12, 11, 11, 10, 10, 10, 7, 8, 8, 8, 8, 6, 5, 8, 8, 9, 10, 11, 11, 8,
    9, 9, 10, 11, 11, 12, 14, 13, 14, 14, 14, 14, 10, 10, 12, 12, 10, 11, 11, 11, 11, 8, 7,
    8, 8, 9, 9, 10, 10, 10, 11, 8, 7, 8, 7, 7, 8, 7, 7, 7, 5, 5, 5, 5, 5, 4, 3, 4, 5, 4,
    3, 4, 5, 5, 4, 4, 5, 6, 5, 5, 4, 4, 4, 2, 2, 4, 2, 1, 0, -1, 0, 1, 1, 1, 1, 0, 0, -1,
    -1, -1, -1, -2, -2, -2, -2, -3, -3, -3, -3, -4, -4, -4, -4, -6, -10, -8, -7, -6, -7, -6,
    -6, -9, -8, -8, -8, -9, -9, -9, -10, -10, -10, -10, -10, -9, -9, -9, -9, -11, -9, -8,
    -10, -10, -10, -11, -11, -10, -10, -11, -11, -12, -12, -14, -12, -10, -10, -9, -10,
    -10, -10, -9, -9, -9, -8, -9, -9, -11, -9, -9, -9, -9, -8, -10, -10, -11, -10, -9, -8,
    -7, -7, -9, -8, -7, -6, -7, -8, -8, -8, -8, -9, -9, -9, -10, -12, -11, -13, -14, -12,
    -13, -12, -12, -11, -12, -13, -11, -11, -13, -14, -14, -14, -14, -13, -12, -13, -12,
    -12, -12, -13, -14, -15, -15, -15, -16, -16, -17, -16, -18, -17, -16, -16, -15, -14,
    -14, -15, -15, -14, -15, -13, -12, -13, -13, -13, -11, -11, -12, -11, -11, -10, -11, -9,
    -9, -9, -9, -11, -11, -9, -11, -8, -9, -7, -10, -11, -10, -12, -13, -15, -15, -16, -16,
    -17, -17, -17, -17, -20, -20, -19, -19, -19, -19, -18, -18, -19, -18, -16, -16, -16,
    -16, -17, -19, -20, -21, -21, -22, -21, -20, -20, -20, -21, -21, -21, -20, -21, -22,
    -23, -23, -25, -24, -24, -25, -24, -24, -25, -24, -25, -24, -22, -21, -21, -23, -24,
    -23, -23, -23, -24, -25, -23, -22, -21, -23, -23, -24, -23, -24, -24, -25, -25, -23,
    -23, -23, -24, -24, -24, -23, -23, -24, -26, -27, -27, -27, -28, -30, -28, -29, -28,
    -27, -27, -27, -28, -29, -31, -31, -31, -31, -32, -31, -30, -30, -31, -33, -34, -35,
    -33, -35, -33, -33, -32, -30, -30, -29, -30, -29, -28, -29, -30, -29, -29, -30, -29,
    -29, -29, -29, -29, -30, -31, -29, -29, -29, -29, -28, -29, -30, -30, -30, -30, -31,
    -31, -32, -32, -33, -33, -33, -33, -32, -32, -34, -34, -35, -34, -34, -34, -35, -35,
    -34, -35, -36, -38, -36, -36, -35, -35, -36, -34, -33, -34, -35, -34, -34, -35, -35,
    -35, -36, -36, -36, -36, -35, -35, -36, -36, -37, -37, -37, -36, -37, -36, -38, -37,

```
   -36, -36, -36, -37, -36, -36, -36, -35, -35, -34, -35, -36, -36, -36, -36, -34, -35,
   -33, -35, -35, -34, -34, -34, -35, -35, -35, -35, -35, -36, -35, -35, -34, -34, -35,
   -37, -38, -38, -38, -38, -37, -37, -38, -38, -37, -38, -38, -39, -40, -39, -36, -38,
   -37, -39, -39, -39, -38, -39];
  lk := floor(1/10000000000000*n) + 1;
  s := 1.000031477579242115061532569306l
      - 0.514839401384136746230446407694400*(1/1000000)*ln(kk);
  s2 := subs(kk = n, s);
  ss := s2*0.9999999999^ll[lk];
  abs(-n*W(-1, -exp(1.0)/n))*ss - piofn
end :
###############################################################
```

Tables of primes :
http://plouffe.fr/NEW/list_primes_pi_of_n_100000000000.txt
http://plouffe.fr/NEW/list_primes_10000000000.txt

Example :
g(1327460000000000,39285023244530) = 49668015014179465.522289485977202
The real value of $p_{1327460000000000} = 49668015014179453$

## Maple program for the calculation of $p_n$, $n = 1.356 \times 10^{16} \dots 10^{24}$

```
###############################################################
F:=proc(n, piofn)
local cn, z, pola, polb;
pola:=.180317882977538680255907226022558343254e-12*x^8-.\
3206852936427839673078154416271278702381e-10*x^7+.\
2521168696363102117361245200766645862014e-8*x^6-.\
1148245660104216214093938301036666192296e-6*x^5+.\
3329760033724798728321163791428487967963e-5*x^4-.\
6343395542494949120689514623217102223176e-4*x^3+.\
7851801857533638277814251770195187519581e-3*x^2-.\
5912431390595948785237451806703967872e-2*x+.\
2187329700777127284427407021768653056673e-1:
polb:=.894992605796963772977753847317326l408730e-12*x^8-.\
1611795950806416304161491053953385128968e-9*x^7+.\
1287542319981049729985262110114907850700e-7*x^6-.\
5988056104228871471776180438025688273194e-6*x^5+.\
1786915791025107343702497983252030617773e-4*x^4-.\
3548565854556946509095877029495597659212e-3*x^3+.\
4690427023808579602996344427037051784331e-2*x^2-.\
398063624996380649025476791478220583827e-1*x+.196997473781\
2788247127674632264178712585:
  z := evalf(log10(n));
  cn := subs(x = z, pola) + subs(x = z, polb)*ln(n);
  evalf(-cn*n*W(-1, -exp(1.0)/n)) - piofn
end:
###############################################################
```

## Table of values of F(n) versus $p_n$

| n | F(n) | $p_n$ |
| --- | --- | --- |
| $10^{16}$ | 394906913903735328.99999995710593 | 394906913903735329 |
| $10^{17}$ | 4185296581467695668.9998280338750 | 4185296581467695669 |
| $10^{18}$ | 44211790234832169331.000076399063 | 44211790234832169331 |
| $10^{19}$ | 465675465116607065549.00000499731 | 465675465116607065549 |
| $10^{20}$ | 4892055594575155744537.0000098572 | 4892055594575155744537 |
| $10^{21}$ | 5127109149801640347 1852.999978699 | 5127109149801640347 1853 |
| $10^{22}$ | 53619387074416211 8627429.00001989 | 53619387074416211 8627429 |
| $10^{23}$ | 559656446798698064 3073682.9999696 | 559656446798698064 3073683 |
| $10^{24}$ | 583100399948365840 70534263.000118 | 583100399948365840 70534263 |

# Table of values of n, $\pi(n)$, $p_n$, G(n)

| $n 10^{15}$ | $\pi(n)$ | $p_n$ | G(n) | Gram Inverted | G(n) $\bar{X} = 2.37605e+07$ | Gap Gram inverted $\bar{X} = 7.62569e+07$ |
|---|---|---|---|---|---|---|
| 1 | 3204941750802 | 3475385758524527 | 3475385752465280 | 3475385760290722 | 6059247 | 1766195 |
| 2 | 6270424651315 | 7093600525704677 | 7093600531547406 | 7093600514882155 | 5842729 | 10822522 |
| 3 | 9287441600280 | 10765662794071351 | 10765662776140778 | 10765662798101237 | 17930573 | 4029886 |
| 4 | 12273824155491 | 14472680634646931 | 14472680642211900 | 14472680659410410 | 7564969 | 24763479 |
| 5 | 15237833654620 | 18205684894350047 | 18205684890027179 | 18205684845589213 | 4322868 | 48760834 |
| 6 | 18184255291570 | 21959393830706447 | 21959393831829265 | 21959393831263666 | 1122818 | 557219 |
| 7 | 21116208911023 | 25730318403586483 | 25730318401099988 | 25730318388727176 | 2496495 | 14859307 |
| 8 | 24035890368161 | 29515978892552597 | 29515978901069447 | 29515978942077307 | 8516850 | 49524710 |
| 9 | 26944926466221 | 33314521777674083 | 33314521779133363 | 33314521740381476 | 1459280 | 37292607 |
| 10 | 29844570422669 | 37124508045065437 | 37124507999149021 | 37124508056355511 | 45916416 | 11290074 |
| 11 | 32735816605908 | 40944788655376237 | 40944788664190013 | 40944788642442631 | 8813776 | 12933606 |
| 12 | 35619471693548 | 44774424266565143 | 44774424274288359 | 44774424246530936 | 7723216 | 20034207 |
| 13 | 38496205973965 | 48612632821248317 | 48612632846598877 | 48612632816905836 | 25350560 | 4342481 |
| 14 | 41366582391891 | 52458753029241283 | 52458753010788072 | 52458753051854838 | 18453211 | 22613555 |
| 15 | 44231080178273 | 56312218341118283 | 56312218348943058 | 56312218418328247 | 7824775 | 77209964 |
| 16 | 47090114439072 | 60172538090123567 | 60172538133649809 | 60172538130369336 | 43526242 | 40245769 |
| 17 | 49944045778207 | 64039282905020807 | 64039282881733481 | 64039282901777277 | 23287326 | 3243530 |
| 18 | 52793190012734 | 67912074089826233 | 67912074089530037 | 67912074062806200 | 296196 | 27020033 |
| 19 | 55637829945151 | 71790575058422851 | 71790575056044677 | 71790575102903791 | 2378174 | 44480940 |
| 20 | 58478215681891 | 75674484987354031 | 75674484995999520 | 75674484998897799 | 8645489 | 11543768 |
| 21 | 61314571044765 | 79563532882638499 | 79563532826010748 | 79563532850159447 | 29570610 | 1727751 |
| 22 | 64147099298639 | 83457473636497967 | 83457473633825154 | 83457473716619734 | 2672813 | 80121767 |
| 23 | 66975984145551 | 87356084739486881 | 87356084732880212 | 87356084781839841 | 6606669 | 42352960 |
| 24 | 69801392791572 | 91259162764140311 | 91259162732170585 | 91259162746099932 | 31969726 | 18040379 |
| 25 | 72623478149504 | 95166521200910351 | 95166521209928786 | 95166521244817945 | 9018435 | 43907594 |
| 26 | 75442380316713 | 99077988824381269 | 99077988830946491 | 99077988840566937 | 6565222 | 16185668 |
| 27 | 78258228083239 | 102993407452131551 | 102993407407852376 | 102993407298472241 | 44279175 | 153659310 |
| 28 | 81071142895913 | 106912630241974061 | 106912630255488703 | 106912630117470118 | 13514642 | 124503943 |
| 29 | 83881233426790 | 110835521231391421 | 110835521272734267 | 110835521271667139 | 41342846 | 40275718 |
| 30 | 86688602810119 | 114761954175793079 | 114761954201984433 | 114761954125523221 | 43638259 | 50269858 |
| 31 | 89493347331727 | 118691810606522897 | 118691810578220391 | 118691810493848076 | 28302506 | 112674821 |
| 32 | 92295556538011 | 122624979771448267 | 122624979772596434 | 122624979823224005 | 1148167 | 51775738 |
| 33 | 95095312182517 | 126561358313013181 | 126561358306791545 | 126561358475859144 | 6221636 | 162845963 |
| 34 | 97892695611204 | 130500849055080079 | 130500849011631709 | 130500849100333580 | 43448370 | 45253501 |
| 35 | 100687778906831 | 134443359677104443 | 134443359776440089 | 134443359790334740 | 10863821 | 109341646 |
| 36 | 103480631416721 | 138388805003019359 | 138388804996907016 | 138388805023564749 | 6112343 | 20545390 |
| 37 | 106271318433884 | 142337102440574897 | 142337102426378912 | 142337102363652797 | 14195985 | 76922100 |
| 38 | 109059901535155 | 146288175048719531 | 146288175029401561 | 146288174931574262 | 19317970 | 117145269 |
| 39 | 111846440164164 | 150241949628632893 | 150241949625523221 | 150241949625523221 | 10535860 | 2180107 |
| 40 | 114630988904000 | 154198357111745083 | 154198357051757170 | 154198357098815921 | 59987913 | 12929162 |
| 41 | 117413599364789 | 158157331435725499 | 158157331422444844 | 158157331469048028 | 13280655 | 33322529 |
| 42 | 120194323133703 | 162118810111632083 | 162118810106649332 | 162118810073330690 | 4982751 | 38301393 |
| 43 | 122973207007771 | 166082733219061063 | 166082733200372932 | 166082733233795825 | 18688131 | 14734762 |
| 44 | 125750296138286 | 170049044111273447 | 170049044050074229 | 170049044050074229 | 21597086 | 61199418 |
| 45 | 128525633848847 | 174017688270699971 | 174017688308542721 | 174017688209805173 | 37842750 | 60894798 |
| 46 | 131299259981906 | 177988613752767869 | 177988613735417195 | 177988613816002938 | 17350674 | 63235069 |
| 47 | 134071214963486 | 181961771195870689 | 181961771178770425 | 181961771229440096 | 17100264 | 33569407 |
| 48 | 136841535130789 | 185937113010539323 | 185937113149378902 | 185937112924455641 | 6248579 | 86083682 |
| 49 | 139610257999130 | 189914593393409557 | 189914593387670038 | 189914593356784763 | 5739519 | 36624794 |
| 50 | 142377417196364 | 193894168896897487 | 193894168924671722 | 193894168842184136 | 27774235 | 54713351 |
| 51 | 145143045599692 | 197875797467973511 | 197875797402159407 | 197875797444772385 | 65814104 | 23201126 |
| 52 | 147907174371027 | 201859438775606323 | 201859438782059357 | 201859438874132705 | 6453034 | 98526382 |
| 53 | 150669836017291 | 205845054246666509 | 205845054239904730 | 205845054390334700 | 39256662 | 143668191 |
| 54 | 153431057455345 | 209832606663909601 | 209832606641202960 | 209832606716127913 | 22706641 | 52218312 |
| 55 | 156190867055604 | 213822059884629769 | 213822059911372304 | 213822059955642707 | 26742535 | 71012938 |
| 56 | 158949293526663 | 217813379543515117 | 217813379570306253 | 217813379519006579 | 26791136 | 24508538 |
| 57 | 161706360526093 | 221806532073166357 | 221806532003780790 | 221806532052347536 | 69385567 | 20818821 |
| 58 | 164462095231054 | 225801485307500567 | 225801485379687884 | 225801485372711222 | 2258919 | 65211155 |
| 59 | 167216521960016 | 229798208333470963 | 229798208351343272 | 229798208407471471 | 17872309 | 74000508 |
| 60 | 169969662554551 | 233796671062467577 | 233796671042312318 | 233796671137843085 | 20155259 | 75375508 |
| 61 | 172721540727639 | 237796844582199317 | 237796844551866864 | 237796844546171680 | 30332453 | 36027637 |
| 62 | 175472177511800 | 241798700605693853 | 241798700580252746 | 241798700646530936 | 25436427 | 39019632 |
| 63 | 178221594869615 | 245802212004481091 | 245802212011032242 | 245802212039361728 | 6551151 | 34880637 |
| 64 | 180969812069916 | 249807352707175657 | 249807352720584206 | 249807352666888216 | 13408549 | 40287441 |
| 65 | 183716850192783 | 253814097057270587 | 253814097060656910 | 253814096974097674 | 3386323 | 83172913 |
| 66 | 186462726814356 | 257822420453703037 | 257822420416882932 | 257822420270061405 | 36820105 | 183641632 |
| 67 | 189207462479325 | 261832298787781253 | 261832298712411022 | 261832298746471022 | 15186608 | 166095511 |
| 68 | 191951073132231 | 265843709118454979 | 265843709097623437 | 265843708773837273 | 20831542 | 344617706 |
| 69 | 194693578185957 | 269856628338594107 | 269856628301551602 | 269856628210471997 | 37042505 | 128122110 |
| 70 | 197434994078331 | 273871034935338403 | 273871034964904223 | 273871035032220227 | 29565820 | 96881824 |
| 71 | 200175335630483 | 277886907995191811 | 277886908023169839 | 277886907974701622 | 27978028 | 20490189 |
| 72 | 202914620525448 | 281904226487172271 | 281904226482872721 | 281904226372807541 | 4145361 | 114364729 |
| 73 | 205652862425306 | 285922970160712259 | 285922970138620930 | 285922970135722026 | 22091329 | 24990233 |
| 74 | 208390079110978 | 289943119715176507 | 289943119700984084 | 289943119723311955 | 14192423 | 8135448 |
| 75 | 211126283162243 | 293964656108882903 | 293964656105170428 | 293964656123793851 | 3712475 | 14910948 |
| 76 | 213861489506392 | 297987560649158759 | 297987560638374821 | 297987560832592311 | 10783938 | 183433431 |
| 77 | 216595711439565 | 302011815896566267 | 302011815820514117 | 302011815832310807 | 76052150 | 64255460 |
| 78 | 219328963332630 | 306037403464228481 | 306037403503476006 | 306037403573745085 | 39247525 | 109516604 |
| 79 | 222061256928013 | 310064306872174139 | 310064306890230891 | 310064306957868189 | 18056752 | 85694050 |
| 80 | 224792606318600 | 314092509321252353 | 314092509284502245 | 314092509318729703 | 36750108 | 2522650 |
| 81 | 227523023099978 | 318121994310899963 | 318121994311396755 | 318121994407209583 | 496792 | 96309620 |
| 82 | 230252520816828 | 322152746418376529 | 322152746409788075 | 322152746375574649 | 8588454 | 42801880 |
| 83 | 232981109132553 | 326184749734071211 | 326184749734869492 | 326184749762788657 | 798281 | 28717446 |
| 84 | 235708800471211 | 330217989477361159 | 330217989421989722 | 330217989480530531 | 55371437 | 3169372 |
| 85 | 238435607431737 | 334252450837522181 | 334252450875892179 | 334252450799878593 | 38369998 | 37643588 |
| 86 | 241161539806582 | 338288119306185553 | 338288119281527478 | 338288119338621616 | 24658075 | 32436063 |
| 87 | 243886608249438 | 342324980816825683 | 342324980787212977 | 342324981049160236 | 29612706 | 232334553 |
| 88 | 246610822221359 | 346363022043182897 | 346363022117594024 | 346363022206964776 | 74411127 | 163781879 |
| 89 | 249334194029546 | 350402229198161719 | 350402229275339670 | 350402229399557886 | 77177951 | 201396167 |
| 90 | 252056733453928 | 354442589428511471 | 354442589461661193 | 354442589515992467 | 33149722 | 87480996 |
| 91 | 254778448725151 | 358484089738515107 | 358484089718761636 | 358484089736797400 | 19753471 | 1717707 |
| 92 | 257499349768637 | 362526717392633107 | 362526717434878969 | 362526717524365342 | 42245862 | 131732235 |
| 93 | 260219446617109 | 366570460630428997 | 366570460634193232 | 366570460613758540 | 3764235 | 16670457 |
| 94 | 262938747499423 | 370615307040895867 | 370615307025898287 | 370615307003910246 | 14997580 | 36985621 |
| 95 | 265657263117161 | 374661245056625183 | 374661245046415361 | 374661244949200614 | 10209822 | 107424569 |
| 96 | 268375000740770 | 378708263125578223 | 378708263087697669 | 378708262951387398 | 37880554 | 174190825 |
| 97 | 271091969073361 | 382756349961937363 | 382756350004680354 | 382756349751872976 | 42742991 | 210064387 |
| 98 | 273808176380030 | 386805494468242607 | 386805494541276472 | 386805494324290341 | 73033865 | 143952266 |
| 99 | 276523631752529 | 390855686180411407 | 390855686125149357 | 390855685867391783 | 55262050 | 313019624 |
| 100 | 279238341033925 | 394906913903735329 | 394906913901321500 | 394906913798224974 | 2413829 | 105510355 |
| 101 | 281952314626716 | 398959167795806791 | 398959167769420669 | 398959167745582098 | 26386122 | 50224693 |
| 102 | 284665559556332 | 403012437560594987 | 403012437526461143 | 403012437543708464 | 34133844 | 16886523 |
| 103 | 287378081126626 | 407066713212243179 | 407066713209546238 | 407066713226257928 | 2696941 | 14014749 |
| 104 | 290089890632238 | 411121985053339019 | 411121985043237519 | 411121985020483115 | 10101500 | 32855904 |
| 105 | 292800991715569 | 415178243427089041 | 415178243406068862 | 415178243341649133 | 21020179 | 85439908 |
| 106 | 295511394070886 | 419235478872659203 | 419235478893349398 | 419235478787660177 | 20690195 | 84999026 |
| 107 | 298221102772488 | 423293682191217493 | 423293682155261572 | 423293682133888953 | 35955921 | 57328540 |
| 108 | 300930125986760 | 427352844259123877 | 427352844262971361 | 427352844328199414 | 3847484 | 69075537 |
| 109 | 303638470099198 | 431412956243052547 | 431412956202837881 | 431412956486153933 | 40214666 | 243101386 |
| 110 | 306346140642929 | 435474009565503149 | 435474009643087789 | 435474009886396362 | 77584640 | 320893213 |
| 111 | 309053145121220 | 439535995988680111 | 439535995939289038 | 439535995966203053 | 49391073 | 22477058 |
| 112 | 311759489314579 | 443598906280158487 | 443598906333361303 | 443598906317194216 | 53202816 | 37035729 |
| 113 | 314465179385261 | 447662732746235623 | 447662732741496788 | 447662732681198514 | 4738835 | 65037109 |
| 114 | 317170221634362 | 451727466937048793 | 451727466961134267 | 451727466946264044 | 24085474 | 9215251 |
| 115 | 319874623177404 | 455793101280411463 | 455793101292794303 | 455793101142809342 | 12382840 | 137602121 |
| 116 | 322578388623503 | 459859627460138027 | 459859627451525756 | 459859627439908299 | 8612271 | 20229728 |
| 117 | 325281523355857 | 463927038138601217 | 463927038105943453 | 463927038141703196 | 32657764 | 3101979 |
| 118 | 327984033568074 | 467995325619454117 | 467995325646787440 | 467995325693804424 | 27433354 | 64486307 |
| 119 | 330685925327709 | 472064482769644943 | 472064482733573233 | 472064482630623649 | 36071710 | 139021294 |
| 120 | 333387204489157 | 476134501830546337 | 476134501867733948 | 476134501670779509 | 37187611 | 159766828 |
| 121 | 336087875323188 | 480205375878386027 | 480205375860189505 | 480205375615331145 | 18196522 | 263054882 |
| 122 | 338787944139611 | 484277097367553351 | 484277097394075116 | 484277097394075116 | 19820714 | 26521765 |
| 123 | 341487414778273 | 488349660001442959 | 488349660004218353 | 488349660052757448 | 2775394 | 51314489 |
| 124 | 344186293058920 | 492423056707800949 | 492423056699489850 | 492423056750244799 | 8311099 | 42443850 |
| 125 | 346884583805017 | 496497280793610557 | 496497280802378550 | 496497280755786896 | 8767993 | 37823661 |
| 126 | 349582292881340 | 500572325543785867 | 500572325541126167 | 500572325446366597 | 2659700 | 97419270 |
| 127 | 352279423771442 | 504648184370381627 | 504648184297048976 | 504648184304134090 | 73332651 | 66247537 |
| 128 | 354975982263335 | 508724850954477793 | 508724850955762227 | 508724850913921935 | 1284434 | 40555858 |
| 129 | 357671973817060 | 512802318988638269 | 512802318948722249 | 512802318960837771 | 39916020 | 27800498 |
| 130 | 360367400804331 | 516880582141749971 | 516880582162249221 | 516880582227931689 | 20499250 | 86181718 |
| 131 | 363062269659721 | 520959634421249321 | 520959634368944099 | 520959634593935391 | 52305222 | 172686070 |
| 132 | 365756583868551 | 525039469767348079 | 525039469802292915 | 525039470031070401 | 34944836 | 263722322 |
| 133 | 368450348555798 | 529120082458008373 | 529120082475909098 | 529120082602922713 | 17900725 | 144914340 |
| 134 | 371143567919892 | 533201466236078989 | 533201466278506017 | 533201466462381357 | 42427028 | 226302368 |
| 135 | 373836245057725 | 537283615564355927 | 537283615609092135 | 537283615849638545 | 44736208 | 285282618 |

# Table of known programs, formulas for the calculation of primes.

| Author (s) | Year | Comment | Efficiency | Number of terms calculated |
|---|---|---|---|---|
| Erathosthène | -276 to -194 | Sieve | Practical | Calculable infinity |
| Mersenne | 1536 | Primes of the form: $2^p$-1. | Practical and exact | 51 |
| Fermat | 1640 | Little theorem of Fermat | Produces weak probable primes | Calculable infinity |
| Euler | 1772 | Second degree polynomial | Practical | 40 |
| Mills | 1947 | Double exponential | Practical | Les than 10 terms |
| Wright | 1951 | Super exponential | Practical | Less than 5 terms |
| Wilson | Circa 1780 | Formula with p ! | Theoretical | Very few |
| Jones, Sato, Wada, Wiens | 1976 | Polynôme de degré 25 à 26 variables | Theoretical | Very few |
| John H. Conway | 1987 | FRACTRAN | Theoretical | Very few |
| Rowland | 2008 | Récurrence | Theoretical | Very few |
| Dress, Landreau | 2010 | 6th degree polynomial | Practical | 58 |
| Benoit Perichon (et al). | 2010 | 26 primes in arithmetic progression | Practical | 26 |
| Kim Walisch | 2020 | Primesieve and Primecount | Practical | Calculable infinity |

# Bibliography (not sorted).

[1] Encyclopedia of Integer Sequences, N.J.A. Sloane, Simon Plouffe, Academic Press , San Diego 1995.
[2] Mills, W. H. (1947), *A prime-representing function, Bulletin of the American Mathematical Society 53 (6): 604, doi:10.1090/S0002-9904-1947-08849-2.*
[3] *E*. M. Wright (1951). *A prime-representing function. American Mathematical Monthly. 58 (9): 616–618. doi:10.2307/2306356. JSTOR 2306356.*
[4] The OEIS, Online Encyclopedia of Integer Sequences, sequences :
 sequences  A051021, A051254, A016104  and A323176, A006988, A006880.
[5] Wikipedia : formulas for primes (effective and non-effective formulas).
https://en.wikipedia.org/wiki/Formula_for_primes
[6] Baillie Robert, The Wright's fourth prime:
 https://arxiv.org/pdf/1705.09741.pdf
[7] Wikipedia : Le recuit simulé : https://fr.wikipedia.org/wiki/Recuit_simul%C3%A9
[8] Wikipedia : Simulated Annealing : https://en.wikipedia.org/wiki/Simulated_annealing
[9] László Tóth, A Variation on Mills-Like Prime-Representing Functions, ArXiv :
https://arxiv.org/pdf/1801.08014.pdf
[10] Makoto Kamada Prime numbers of the form 7, 73, https://stdkmd.net/nrr/7/73333.htm#prime
[11] Plouffe, Simon : Pi, the primes and the Lambert W function, conference in July 2019, Montréal at the ACA 2019 (ETS). https://vixra.org/abs/1907.0108
[12] Gram, J. P. "Undersøgelser angaaende Maengden af Primtal under en given Graeense." *K. Videnskab. Selsk. Skr.* **2**, 183-308, 1884.
[13] Kim Walisch, primecount and primesieve, fastest program to compute primes.
https://github.com/kimwalisch/primecount
[14] Visser, Matt  : Primes and the Lambert W function : https://arxiv.org/abs/1311.2324
[15] Berndt, Bruce, Ramanujan Notebooks IV, page 124.
[16 ] Dusart, Pierre, Autour de la fonction qui compte le nombre de nombres premiers, thèse de Doctorat 1998.
[17] Plouffe, Simon, List of primes computed by the primecount program :
http://plouffe.fr/NEW/list_primes_pi_of_n_100000000000.txt
http://plouffe.fr/NEW/list_primes_10000000000.txt
[18] Kahane, Jean-Pierre, le nombre cet inconnu. : http://ww3.ac-poitiers.fr/math/prof/resso/kah/conf.pdf
[19] Skewe's Number on wikipedia https://en.wikipedia.org/wiki/Skewes%27s_number
[20] Salvy, Bruno, Fast computation of some asymptotic functional inverses, J. Symbolic Comput.17(1994), 227–236
[21] Mendès-France, Michel, Tannenbaum Les nombres premiers, entre l'ordre et le chaos.
[22] Plouffe, Simon , Pi, the primes and the Lambert W function. Conférence ACA 2019, Montréal ETS.
https://vixra.org/abs/1907.0108
[23] Plouffe, Simon , Nombre premiers en progression géométrique : http://plouffe.fr/NEW/
[24] Weisstein, Eric W. "Gram Series." Mathworld at https://mathworld.wolfram.com/GramSeries.html